\documentclass{amsart}
\usepackage{amssymb}
\usepackage{amsmath}
\usepackage{amsthm}
\usepackage{longtable}
\usepackage{graphicx}
\usepackage{float}
\usepackage{longtable}
\newcommand{\pn}{\par\noindent}
\newcommand{\pmn}{\par\medskip\noindent}
\newtheorem{theor}{Theorem}[section]
\theoremstyle{definition} 
\newtheorem{ex}{Example}[section] \theoremstyle{remark}
\newtheorem{rem}{Remark}[section]
\begin{document}
\title{On the geometry of real polynomials of the 7-th degree with
simple complex critical points}
\author{Elizaveta Khokhlova and Yury Kochetkov}
\date{}
\email{elankhokhlova@edu.hse.ru,yukochetkov@hse.ru,yuyukochetkov@gmail.com}

\begin{abstract} We consider real monic polynomials with simple complex
critical points. Critical values of such polynomial $p$ can be
considered as vertices of the uniquely defined polygon, symmetric
with respect to the $x$-axis. The inverse image of this polygon is
a tree-like system of ovals that touch each other externally in
critical points. This system of ovals will be called the
geometrical code of $p$. Let $a_1+i\cdot b_1,\ldots,a_n+i\cdot
b_n$, $a_1<a_2<\ldots<a_n$ and $b_i>0$, $i=1,\ldots,n$, be the set
of critical points of our polynomial $p$ in the upper halfplane,
and $c_1+i\cdot d_1,\ldots,c_n+i\cdot d_n$  be their images. The
order of numbers $c_k$ and the list of signs of numbers $d_k$ will
be called the algebraic code of $p$. In our experimental work we
study geometrical and algebraic codes of the 7-th degree
polynomials and the relations between them.
\end{abstract}

\maketitle

\section{Introduction}
\pn We want to correspond a figure to a polynomial
$p\in\mathbb{C}[z]$. It can be done in the following way: we draw
an oval $O$ that passes through all critical values of $p$ and
consider the inverse image $T=p^{-1}(O)$. $T$ is a tree-like
structure of ovals that externally touch each  in critical points.
There are two well-known cases. The first case is when $p$ has
exactly two finite critical values $\alpha$ and $\beta$. Then $O$
is just a segment $[\alpha,\beta]$ and $T$ is a plane tree. Such
polynomials are called generalized Chebyshev polynomials or Shabat
polynomials. They are studied in the scope of the dessins
d'enfants theory \cite{LZ}.
\begin{ex} Let $p=4z^3-3z^4$. Here $0$ and $1$ are critical points and
$0$ and $1$ --- critical values. $O=[0,1]$,
\[\begin{picture}(100,30) \put(15,4){\tiny $-\frac 13\,-\frac{\sqrt{2}\,i}{3}$}
\put(15,24){\tiny $-\frac 13\,+\frac{\sqrt{2}\,i}{3}$}
\multiput(55,5)(0,20){2}{\circle*{2}} \put(55,5){\line(1,1){10}}
\put(55,25){\line(1,-1){10}} \put(65,15){\line(1,0){25}}
\multiput(65,15)(25,0){2}{\circle*{2}} \put(80,15){\circle*{2}}
\put(65,8){\tiny $0$} \put(79,8){\tiny $1$} \put(95,13){\tiny
$\frac 43$} \put(0,13){\tiny $T=$} \end{picture}\]
\begin{center}{Figure 1}\end{center} \end{ex}
\pn The second case is when $p$ has exactly three finite critical
values $\alpha$, $\beta$ and $\gamma$. Then $O$ is the triangle
$\alpha\beta\gamma$ and $T$ is a cactus. Such polynomials are
called Zolotarev polynomials. They are also investigated in the
scope of the dessins d'enfants theory \cite{LZ},\cite{YK}.
\begin{ex} Let
$$p=\int(x^2+1)\cdot(x-1)\cdot\left(x+\frac 23-\frac{\sqrt{35}}{3}\,i\right)
\,dx.$$ Here $1,\pm i,-2/3+\sqrt{35}\,i/3$ are critical points and
$p\,(1)=p\,(-2/3+\sqrt{35}\,i/3)$. The corresponding cactus is of
the form:
\[\begin{picture}(80,35) \put(0,10){\line(1,0){80}}
\multiput(0,10)(70,-10){2}{\line(1,1){10}}
\put(10,20){\line(1,-1){20}}
\put(30,0){\line(1,1){30}} \put(40,30){\line(1,-1){30}}
\put(40,30){\line(1,0){20}} \end{picture}\]
\begin{center}{Figure 2}\end{center} \end{ex}
\pn If $p$ is a generic polynomial, then the oval $O$ cannot be
defined uniquely and different ovals generate different tree-like
structures $T$. \pmn However, if $p$ is a real polynomial of an
odd degree with simple complex critical points, then the set of
its critical values is symmetric with respect to the $x$-axis.
Thus, there exists the unique polygon with vertices in critical
values:
\[\begin{picture}(230,70) \put(0,40){\vector(1,0){90}}
\put(10,10){\vector(0,1){60}} \put(3,67){\tiny y}
\put(87,33){\tiny x} \multiput(20,30)(0,20){2}{\circle*{2}}
\multiput(40,20)(0,40){2}{\circle*{2}}
\multiput(70,30)(0,20){2}{\circle*{2}} \put(20,5){\small Critical
values}

\put(110,37){$\Rightarrow$}

\put(140,40){\vector(1,0){90}} \put(150,10){\vector(0,1){60}}
\put(143,67){\tiny y} \put(227,33){\tiny x}
\multiput(160,30)(0,20){2}{\circle*{2}}
\multiput(180,20)(0,40){2}{\circle*{2}}
\multiput(210,30)(0,20){2}{\circle*{2}}
\put(160,30){\vector(2,-1){20}} \put(180,20){\vector(3,1){30}}
\put(210,30){\vector(0,1){20}} \put(210,50){\vector(-3,1){30}}
\put(180,60){\vector(-2,-1){20}} \put(160,50){\vector(0,-1){20}}
\put(170,5){\small Polygon}
\end{picture}\]
\begin{center}{Figure 3}\end{center}
The inverse image of this polygon is a tree-like structure of
ovals that touch each other externally and symmetric with respect
the the $x$-axis. \pn In what follows we will consider real
polynomials of odd degree with simple complex critical points.
Let $p$ be such polynomial of degree $2n+1$ and let $a_1\pm
i\,b_1,\ldots,a_n\pm i\, b_n$, $a_1<a_2<\ldots<a_n$, $b_1>0,\ldots,b_n>0$,
be its critical points. Let $p(a_k+b_k\,i)=c_k+i\,d_k$, $k=1,\ldots n$,
and let $c_k\neq c_l$, $k\neq l$, and $d_k\neq 0$, $k=1,\ldots,n$.

\begin{rem} Without loss of generality we can assume that $a_1=0$
and $a_n=1$. \end{rem} \pn The increasing order of numbers $c_i$
defines a permutation $\sigma\in S_n$. This permutation with the
list of signs of numbers $d_i$ constitute the \emph{algebraic
code} of a polynomial $p$.
\begin{ex} Let $p=7\cdot\int((x^2+1)\cdot(x^2-\frac 23\cdot x+\frac
59)\cdot (x^2-2x+2)\,dx$. Here $\pm i,\frac 13\pm \frac 23\cdot i$
and $1\pm i$ are critical points.
$$p(i)\approx 2.72+1.95\,i,\,\, p(1/3+2i/3)\approx 2.13+2.04,\,\,
p(1+i)\approx 0.84+3.51\,i.$$ Thus, $[(3,2,1),[+++]]$ is the
algebraic code of $p$. \end{ex} \pn For polynomials of degree
$2n+1$ there are $2^n\cdot n!$ possible algebraic codes. We
introduce the $n!\times 2^n$ \emph{table of possible codes}. For
$n=7$ the table is of the form:
\[{\scriptsize \begin{tabular}{|c|c|c|c|c|c|c|c|c|} \hline &
[+\,+\,+] & [+\,+\,--] &
[+\,--\,+] & [+\,--\,--] & [--\,+\,+] & [--\,+\,--] &
[--\,--\,+] & [--\,--\,--]\\
\hline (123)
&&&&&&&& \\ \hline (132) &&&&&&&& \\ \hline (213) &&&&&&&& \\
\hline (231) &&&&&&&&\\ \hline (312) &&&&&&&& \\ \hline (321)
&&&&&&&& \\ \hline \end{tabular}}\]
\begin{center}{\small Table 1}\end{center}
\pn Now we will explain how to code the oval structure
$T=p^{-1}(O)$. If $p$ is of degree $2n+1$, then $T$ is a tree-like
structure of $2n+1$ ovals, symmetric with respect to the $x$-axis.
Our coding will take into account only ovals in the upper
half-plane. Each oval contains $2n$ points --- inverse images of
critical values. The counterclockwize order of these points in
each oval is the same as the order in $O$. The central oval of the
structure $T$ is symmetric with respect to the $x$-axis and
contains $n$ points in the upper half-plane and $n$ points in the
lower. We numerate upper points in the central oval from left to
right. At least two points (symmetric) are critical points where
other ovals touch the central. \pmn The central oval has level
zero. The first level ovals have a common point with the central.
The second level ovals have a common point with some first level
oval and so on. The base point of an oval of level $k$, $k>0$, is
its common point with oval of level $k-1$. Numeration of points in
oval begins from the point next to the base in the
counterclockwize order.
\begin{ex}
\[\begin{picture}(230,95)  \put(20,0){\oval(40,60)[lt]}
\put(20,30){\line(1,0){190}} \put(210,0){\oval(40,60)[rt]}
\put(90,0){\small Central oval}
\multiput(40,30)(50,0){4}{\circle*{2}} \put(39,20){\small 1}
\put(89,20){\small 2} \put(139,20){\small 3} \put(189,20){\small
4} \put(40,40){\oval(30,20)[b]} \put(40,60){\oval(30,20)[t]}
\multiput(25,40)(30,0){2}{\line(0,1){20}}
\multiput(25,40)(0,10){3}{\circle*{2}}
\multiput(55,40)(0,10){3}{\circle*{2}} \put(40,70){\circle*{2}}
\put(49,38){\tiny 1} \put(49,48){\tiny 2} \put(49,58){\tiny 3}
\put(39,64){\tiny 4} \put(28,58){\tiny 5} \put(28,48){\tiny 6}
\put(28,38){\tiny 7}

\put(140,40){\oval(30,20)[b]} \put(140,60){\oval(30,20)[t]}
\multiput(125,40)(30,0){2}{\line(0,1){20}}
\multiput(125,40)(0,10){3}{\circle*{2}}
\multiput(155,40)(0,10){3}{\circle*{2}} \put(140,70){\circle*{2}}
\put(149,38){\tiny 1} \put(149,48){\tiny 2} \put(149,58){\tiny 3}
\put(139,64){\tiny 4} \put(128,58){\tiny 5} \put(128,48){\tiny 6}
\put(128,38){\tiny 7}

\put(165,60){\oval(20,20)[l]}
\multiput(165,50)(0,20){2}{\line(1,0){20}}
\put(185,60){\oval(20,20)[r]}
\multiput(165,50)(10,0){3}{\circle*{2}}
\multiput(165,70)(10,0){3}{\circle*{2}} \put(195,60){\circle*{2}}
\put(164,43){\tiny 1} \put(174,43){\tiny 2} \put(184,43){\tiny 3}
\put(184,63){\tiny 5} \put(174,63){\tiny 6} \put(164,63){\tiny 7}
\put(198,58){\tiny 4}

\put(175,75){\oval(10,10)[b]}
\multiput(170,75)(10,0){2}{\line(0,1){12}}
\put(175,87){\oval(10,10)[t]}
\multiput(170,75)(0,6){3}{\circle*{2}}
\multiput(180,75)(0,6){3}{\circle*{2}} \put(175,92){\circle*{2}}

\put(37,47){\scriptsize O} \put(43,45){\tiny 1}
\put(137,47){\scriptsize O} \put(143,45){\tiny 2}
\put(172,56){\scriptsize O} \put(178,54){\tiny 3}
\end{picture}\]
\pmn
\begin{center}{Figure 4}\end{center}
\pn The above structure has code $(1,3(3(6)))$. It must be noted
that this structure cannot be realized as the inverse image of
$O$: the point, symmetric to the point 1 in the central oval and
point 3 in $O_2$ are inverse images of the same point in $O$ (they
are three steps from point 3 in the central oval). But both are
critical points --- impossible.
\end{ex}
\begin{ex}
\[\begin{picture}(320,90)  \put(20,0){\oval(40,60)[lt]}
\put(20,30){\line(1,0){110}} \put(130,0){\oval(40,60)[rt]}
\put(45,0){\small Code $(1,3(3,6))$}
\multiput(30,30)(30,0){4}{\circle*{2}} \put(29,20){\small 1}
\put(59,20){\small 2} \put(89,20){\small 3} \put(119,20){\small 4}
\put(30,40){\oval(10,20)}

\put(90,40){\oval(30,20)[b]} \put(90,60){\oval(30,20)[t]}
\multiput(75,40)(30,0){2}{\line(0,1){20}}
\multiput(75,40)(0,10){3}{\circle*{2}}
\multiput(105,40)(0,10){3}{\circle*{2}} \put(90,70){\circle*{2}}
\put(99,38){\tiny 1} \put(99,48){\tiny 2} \put(99,58){\tiny 3}
\put(88,64){\tiny 4} \put(78,58){\tiny 5} \put(78,48){\tiny 6}
\put(78,38){\tiny 7}

\put(115,60){\oval(20,10)} \put(65,50){\oval(20,10)}

\put(190,0){\oval(40,60)[lt]} \put(190,30){\line(1,0){110}}
\put(300,0){\oval(40,60)[rt]} \put(215,0){\small Code
$(1(4),3(2))$} \multiput(200,30)(30,0){4}{\circle*{2}}
\put(199,20){\small 1} \put(229,20){\small 2} \put(259,20){\small
3} \put(289,20){\small 4}

\put(200,40){\oval(30,20)[b]} \put(200,60){\oval(30,20)[t]}
\multiput(185,40)(30,0){2}{\line(0,1){20}}
\multiput(185,40)(0,10){3}{\circle*{2}}
\multiput(215,40)(0,10){3}{\circle*{2}} \put(200,70){\circle*{2}}
\put(209,38){\tiny 1} \put(209,48){\tiny 2} \put(209,58){\tiny 3}
\put(198,64){\tiny 4} \put(188,58){\tiny 5} \put(188,48){\tiny 6}
\put(188,38){\tiny 7}

\put(260,40){\oval(30,20)[b]} \put(260,60){\oval(30,20)[t]}
\multiput(245,40)(30,0){2}{\line(0,1){20}}
\multiput(245,40)(0,10){3}{\circle*{2}}
\multiput(275,40)(0,10){3}{\circle*{2}} \put(260,70){\circle*{2}}
\put(269,38){\tiny 1} \put(269,48){\tiny 2} \put(269,58){\tiny 3}
\put(258,64){\tiny 4} \put(248,58){\tiny 5} \put(248,48){\tiny 6}
\put(248,38){\tiny 7}

\put(200,80){\oval(10,20)} \put(285,50){\oval(20,10)}
\end{picture}\]
\pmn
\begin{center}{Figure 5}\end{center}
\end{ex}
\pn In what follows we will study algebraic and geometric codes of
polynomials of degree 5 and 7.

\section{Polynomials of degree 5}
\pn Our polynomial $p$ of degree 5 has four critical points and
four critical values. We assume that $\pm a\,i$ and $1\pm b\,i$,
$a>0$, $b>0$, are critical points. Let $c_1+d_1i$, and $c_2+d_2i$
be corresponding critical values. Thus,
$$p=5\int (x^2+a^2)(x^2-2x+1+b^2)\,dx.$$ The first quadrant is the
space of parameters. Curves $d_1=0,d_2=0$ and $c_1-c_2=0$ divide
this space into five domains:
\[\begin{picture}(120,100) \put(0,10){\vector(1,0){100}}
\put(95,2){\tiny a} \put(10,0){\vector(0,1){100}} \put(2,93){\tiny
b} \qbezier(30,10)(50,30)(90,40) \put(94,38){\tiny $d_1=0$}
\qbezier(10,30)(30,50)(40,90) \put(30,93){\tiny $d_2=0$}
\qbezier(20,10)(30,30)(90,60) \qbezier(10,20)(30,30)(60,90)
\put(75,75){\tiny $c_1=c_2$} \put(85,65){$\downarrow$}
\put(60,74){$\leftarrow$} \put(20,70){\scriptsize 1}
\put(40,70){\scriptsize 2} \put(70,20){\scriptsize 5}
\put(70,40){\scriptsize 4} \put(55,55){\scriptsize 3}
\end{picture}\]
\begin{center}{Figure 6}\end{center}
\pmn We have:
\begin{itemize}
\item in the domain 1 $c_1>c_2$, $d_1>0$, $d_2<0$, geometrical
code $(2(2))$; \item in the domain 2 $c_1>c_2$, $d_1>0$, $d_2>0$,
geometrical code $(2(1))$; \item in the domain 3 $c_1<c_2$,
$d_1>0$, $d_2>0$, geometrical code $(1,2)$; \item in the domain 4
$c_1>c_2$, $d_1>0$, $d_2>0$, geometrical code $(1(3))$; \item in the
domain 5 $c_1>c_2$, $d_1<0$, $d_2>0$, geometrical code $(1(2))$.
\end{itemize}
\pn Thus, all possible geometrical codes are realized and the
table of algebraic codes is of the form:
\[{\scriptsize \begin{tabular}{|c|c|c|c|c|} \hline
&[+\,+]&[+\,--]&[--\,+]&[--\,--]\\ \hline (1,2)&\checkmark&0&0&0\\
\hline (2,1)&\checkmark&\checkmark&\checkmark&0\\ \hline
\end{tabular}}\]
\begin{center}{Table 2}\end{center}
\pmn Here symbol "\checkmark" denotes that this code is realized
and $0$ --- not.
\begin{ex} In the picture below the structure $T$ is presented for
$a=2$, $b=1$:
\begin{figure}[H]
    \centering
    \includegraphics[width=0.4\textwidth]{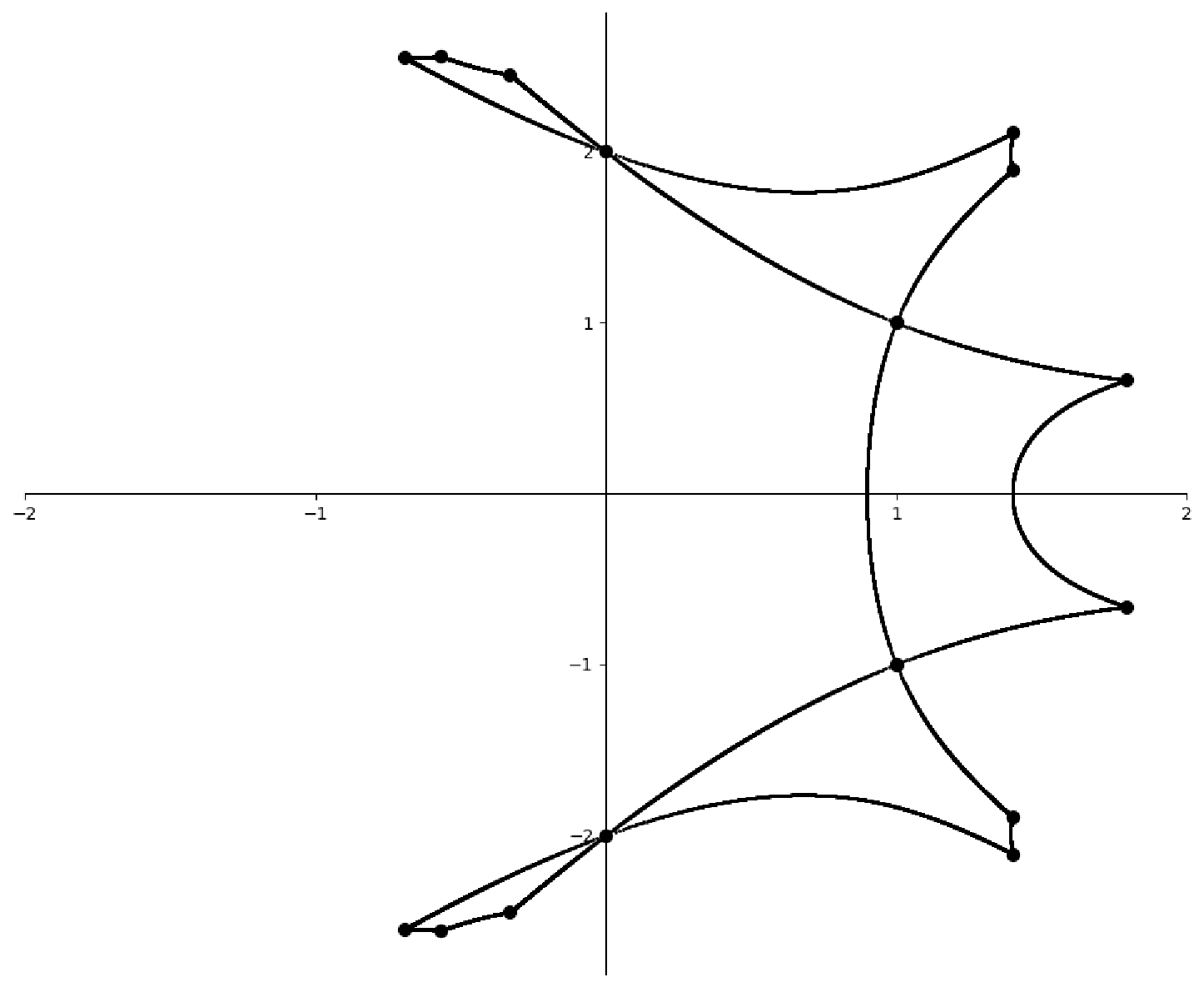}
\end{figure}
\begin{center}{Figure 7}\end{center}
\end{ex}
\begin{rem} Non smooth points in the Figure above are exactly inverse images
of critical values.
\end{rem}
\section{Polynomials of degree 7}
\pn The problem of describing the geometrical structure of
polynomials of  degree 7 is substantially more complicated. After an affine
transformation (translation and scaling), we may assume that the
critical points of the polynomial are of the form $\pm a i$,
$b\pm c i$, $1\pm d i$, where $a>0$, $c>0$, $b\in(0,1)$, and $d>0$.
\pn The corresponding critical values are
 $r_1\pm i q_1$, $r_2\pm i q_2$, and $r_3\pm i q_3$.
\pn Possible algebraic codes are presented in Table 1. \pn
Extensive computations give the following filling of this table:
\[{\scriptsize \begin{tabular}{|c|c|c|c|c|c|c|c|c|} \hline &
[+\,+\,+] & [+\,+\,--] & [+\,--\,+] & [+\,--\,--] & [--\,+\,+] &
[--\,+\,--] &
[--\,--\,+] & [--\,--\,--]\\
\hline (123)&\checkmark&\checkmark&0&0&\checkmark&\checkmark&0&0 \\
\hline (132)&\checkmark&\checkmark&\checkmark&0&\checkmark&
\checkmark&\checkmark&0 \\
\hline (213)&\checkmark&\checkmark&\checkmark&\checkmark&
\checkmark&\checkmark&0&0 \\
\hline (231)&\checkmark&\checkmark&\checkmark&\checkmark&
\checkmark&\checkmark&0&0 \\
\hline (312)&\checkmark&\checkmark&\checkmark&0&\checkmark&
\checkmark&\checkmark&0 \\
\hline (321)&\checkmark&\checkmark&\checkmark&\checkmark&
\checkmark&\checkmark&\checkmark&0 \\
\hline \end{tabular}}\]
\begin{center}{\small Table 3}\end{center}
\pn A quite nontrivial (and AI assisted) reasoning prove the
following statement:
\begin{theor} If $q_2<0$, then either $q_1>0$ and
$r_1>r_2$, or $q_3>0$ and $r_2>r_3$. \end{theor} \pn This
statement explains the form of Table 3. The proof is in Supplement
1. \pmn To understand the geometry of the 7-degree polynomials, we
must study how hypersurfaces $r_1-r_2=0$, $r_1-r_3=0$,
$r_2-r_3=0$, $q_1=0$, $q_2=0$ and $q_3=0$ divide the
four-dimensional parameter space $a,b,c,d$ into domains. We study
this partition using two-dimensional sections of the parameter
space by fixing parameters $a$ and $d$. The section of the
parameter space by the plane $a=d=1$ is presented in the Figure
below.
\[\begin{picture}(250,235)  \put(0,-5){\vector(0,1){240}}
\put(-6,230){\scriptsize y} \qbezier(0,150)(90,170)(100,230)
\put(-5,-2){\vector(1,0){255}} \put(248,-8){\scriptsize x}
\multiput(240,0)(0,10){23}{\line(0,1){6}}
\put(240,-2){\circle*{2}}\put(239,-9){\tiny 1}
\qbezier(240,150)(150,170)(140,230)

\put(0,150){\line(5,-1){240}} \put(240,150){\line(-5,-1){240}}

\qbezier(0,150)(100,70)(100,0) \qbezier(240,150)(140,70)(140,0)

\qbezier(120,126)(60,137)(0,77) \qbezier(120,126)(180,137)(240,77)

\qbezier(0,20)(125,45)(240,40) \qbezier(240,20)(115,45)(0,40)

\put(40,15){\normalsize \textbf{v}} \put(20,30){\normalsize
\textbf{w}} \put(60,0){\small \textbf{x}} \put(120,20){\normalsize
\textbf{y}} \put(170,0){\small \textbf{e}}
\put(220,30){\normalsize \textbf{f}} \put(200,10){\normalsize
\textbf{g}} \put(190,80){\normalsize \textbf{h}}
\put(120,100){\normalsize \textbf{i}} \put(50,80){\normalsize
\textbf{j}}

\put(120,150){\normalsize \textbf{k}} \put(5,90){\normalsize
\textbf{l}} \put(3,130){\normalsize \textbf{m}}
\put(20,110){\normalsize \textbf{n}} \put(50,125){\normalsize
\textbf{o}} \put(20,200){\normalsize \textbf{p}}

\put(233,129){\normalsize \textbf{q}} \put(220,200){\normalsize
\textbf{r}} \put(220,110){\normalsize \textbf{s}}
\put(220,95){\normalsize \textbf{t}} \put(185,125){\normalsize
\textbf{u}}

\put(15,130){\vector (3,-1) {45}} \put(231,132){\vector (-3,-1)
{51}} \put(65,5){\vector (1,1) {32}} \put(173,8){\vector (-1,1)
{30}} \end{picture} \]
\begin{center}{Figure 8}\end{center}

\pn In the table below we demonstrate which algebraic and geometric
codes correspond to each domain.
\[{\scriptsize
\begin{tabular}{|c|c|c|}
\hline

\textbf{Domain} & \textbf{Algebraic code} &
\textbf{Geometric code} \\
\hline e & (321)[+\,+\,--] & (2(2,5)) \\
\hline f & (231)[+\,+\,--] & (1(2,4)) \\
\hline g & (231)[--\,+\,--] & (1(2,3)) \\
\hline h & (231)[+\,+\,+] & 1(5(5)) \\
\hline i & (321)[+\,+\,+] & (2(1,5)) \\
\hline j & (312)[+\,+\,+] & (3(1(1))) \\
\hline k & (123)[+\,+\,+] & (1,2,3) \\
\hline l & (132([+\,+\,+] & (3(1(1))) \\
\hline m & (132)[+\,+\,+] & (1(2),2) \\
\hline n & (123)[+\,+\,+] & (2(1),3) \\
\hline o & (213)[+\,+\,+] & (1(5),3) \\
\hline p & (213)[+\,+\,+] & (2(1),3) \\
\hline q & (231)[+\,+\,+] & (2,3(2)) \\
\hline r & (132)[+\,+\,+] & (1,2(5)) \\
\hline s & (123)[+\,+\,+] & (1,2(5)) \\
\hline t & (213)[+\,+\,+] & (1(4,5)) \\
\hline u & (132)[+\,+\,+] & (1,2(5)) \\
\hline v & (312)[--\,+\,--] & (3(3,4)) \\
\hline w & (312)[--\,+\,+] & (3(2,4)) \\
\hline x & (321)[--\,+\,--] & (2(1,4)) \\
\hline y & (321)[--\,+\,+] & (2(2,4)) \\
\hline \end{tabular}} \]
\begin{center}{Table 4}\end{center}

\section{Geometric codes of the 7th degree polynomials}
\pn Below is the complete list of (possible) geometric codes of the 7th
degree polynomials.
\begin{enumerate}
\item (1,2,3);
\item (1(3),2), (1(4),2), (1(2),3), (1(5),3), (1,2(4)), (1,2(5)), (2(1),3),
(2(2),3), (1,3(1)), (1,3(4)), (2,3(2)), (2,3(3));
\item (1(2,3)), (1(2,4)), (1(3,5)), (1(4,5)), (2(1,4)), (2(1,5)), (2(2,4)),
(2(2,5)), (3(1,2)), (3(1,3)), (3(2,4)), (3(3,4));
\item (1(2(1))), (1(2(2))), (1(3(2))), (1(3(5))), (1(4(1))), (1(4(4))),
(1(5(4))), (1(5(5)));
\item (2(1(3))), (2(1(4))), (2(2(2))), (2(2(3))), (2(4(3))), (2(4(4))),
(2(5(2))), (2(5(3)));
\item (3(1(1))), (3(1(2))), (3(2(2))), (3(2(5))), (3(3(1))), (3(3(4))),
(3(4(4))), (3(4(5))).   \end{enumerate}

\begin{ex} And here are two examples of tree-like structures of ovals.
\[\parbox{6cm}{\includegraphics[scale=0.3]{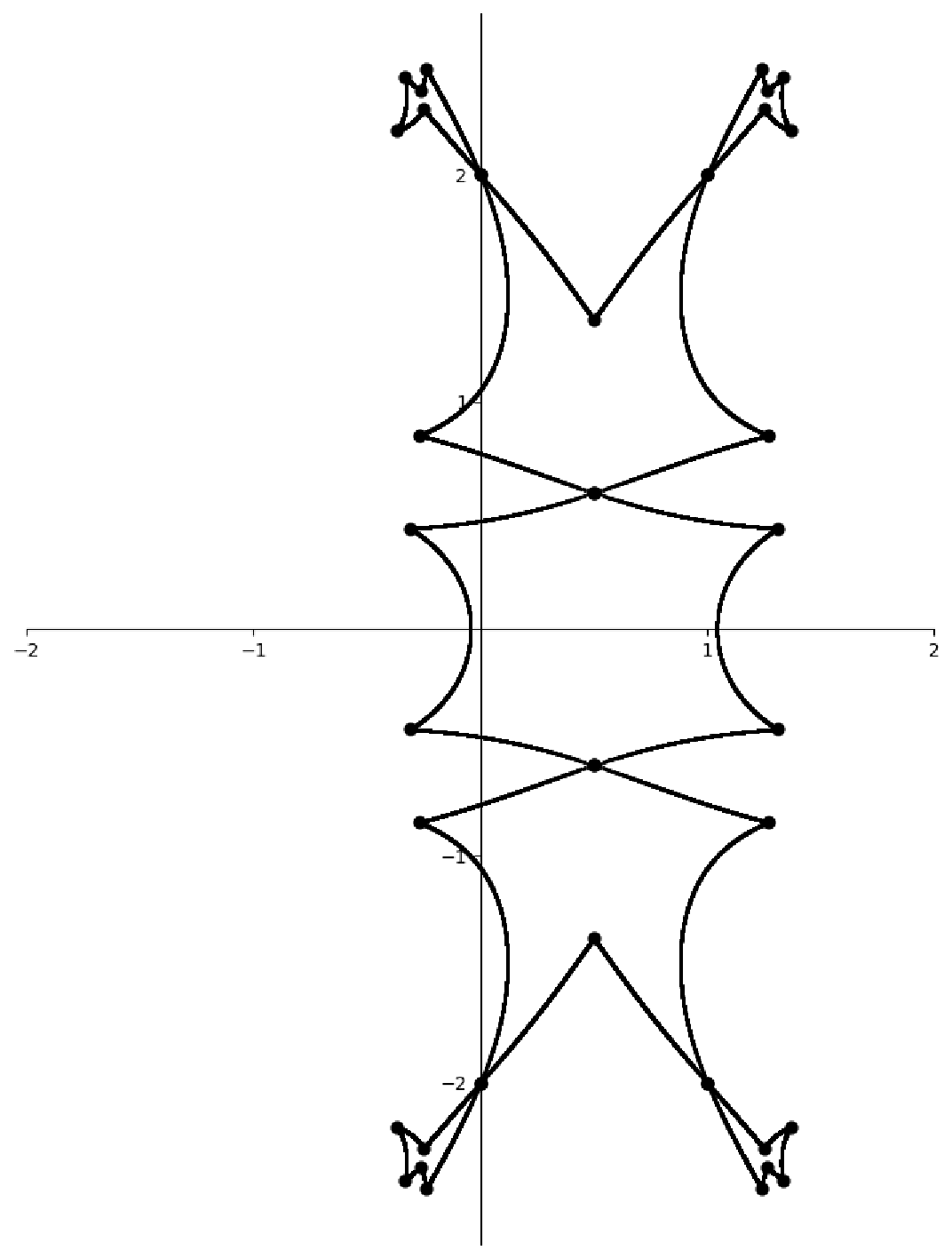}}\quad
\parbox{6cm}{\includegraphics[scale=0.3]{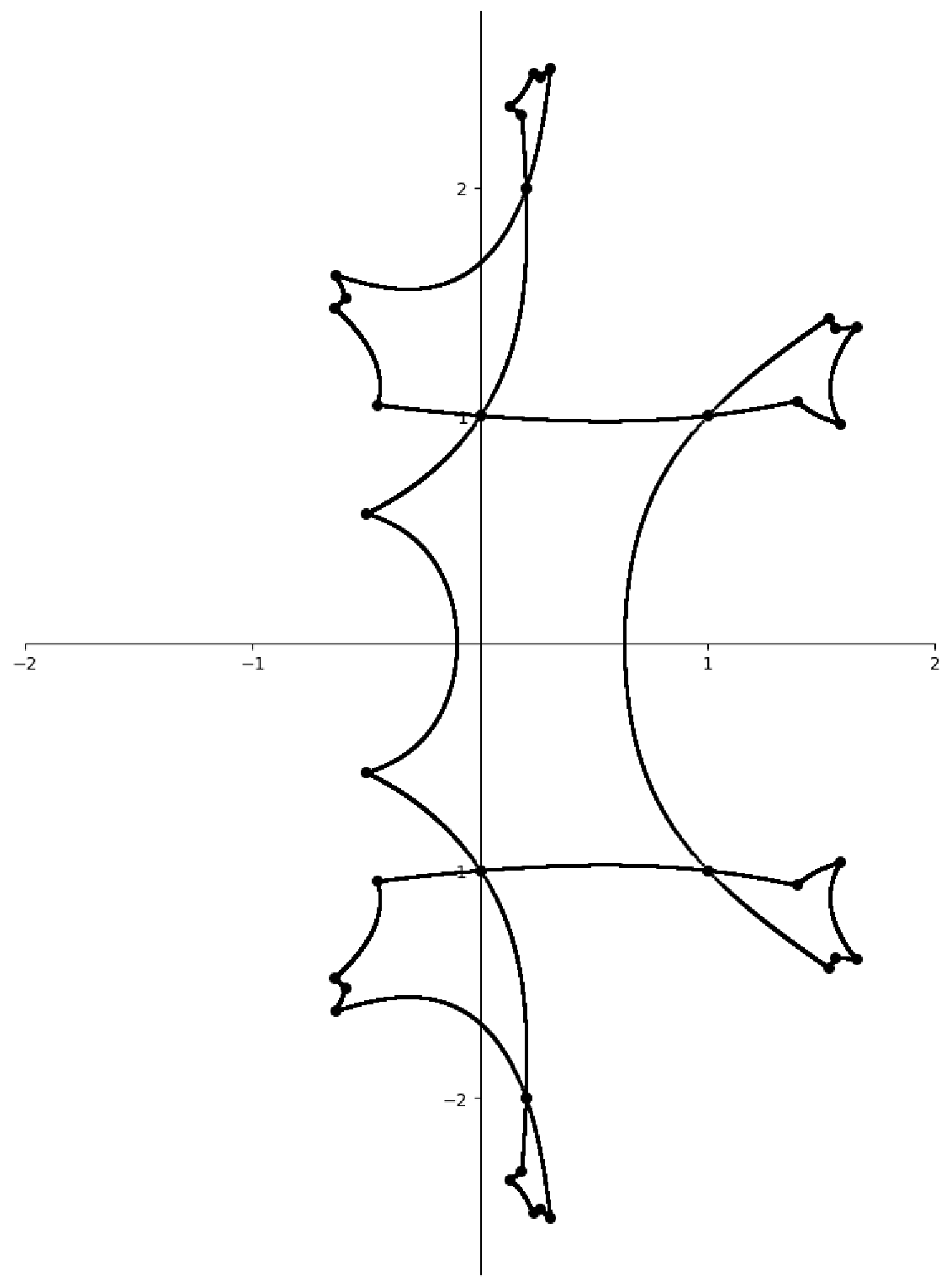}}\]
\pn\hspace{40mm}  $a=d=2,b=0.5,c=0.6$\hspace{27mm}
$a=d=1,b=0.2,c=2$.
\begin{center}{Figure 9}\end{center}
\end{ex}

\pn All possible geometric codes are realized. In the following
table we correspond to each code the values of parameters of its
realization. It must be noted that if some code is realized for
parameters $a_1,b_1,c_1,d_1$, then mirror-symmetric code is
realized for parameters $a_2=d_1,b_2=1-b_1,c_2=c_1,d_2=a_1$.
{\scriptsize \[\begin{tabular}{|c|c|c|c|c|c|} \hline geometric
code&a&b&c&d&
algebraic code\\ \hline (1,2,3)& 0.1&0.1&0.1&0.1&(123)[+\,+\,+]\\
\hline (1(3),2)&02&0.18978&1.934&2&(231)[+\,+\,+]\\ \hline
(1(4),2)&2&0.75&1.9&1.8&(231)[+\,+\,+]\\ \hline
(1(2),3)&2&0.64&1.175&1&(213)[--\,+\,+]\\ \hline
(1(5),3)&2&0.175&1.95&2&(213)[+\,+\,+]\\ \hline
(1,2(4))&0.8&0.9&1.5&0.2&(132)[+\,--\,+]\\ \hline
(1,2(5))&2&0.8&4&2&(132)[+\,+\,+]\\ \hline
(1(2,3))&2&0.8&0.7&2&(231)[--\,+\,--]\\ \hline
(1(2,4))&2&0.9694&0.731&2&(231)[+\,+\,--]\\ \hline
(1(3,5))&3&0.97&3.6&1.2&(312)[+\,--\,+]\\ \hline
(1(4,5))&2&0.980165&1.997&1.9&(312)[+\,+\,+]\\ \hline
(2(1,4))&2&0.36&0.86&2&(321)[--\,+\,+]\\ \hline
(2(1,5))&1&0.5&0.9&1&(321)[+\,+\,+]\\ \hline
(2(2,4))&2&0.5&0.6&2&(321)[--\,+\,--]\\ \hline
(1(2(1)))&0.4&0.1&0.1&1.7&(213)[--\,+\,--]\\ \hline
(1(2(2)))&4&0.56&0.1&1.4&(132)[+\,+\,--]\\ \hline
(1(3(2)))&4&0.015&4.29&1.2&(321)[--\,+\,+]\\ \hline
(1(3(5)))&1.7&0.8&0.7&0.1&(312)[--\,--\,+]\\ \hline
(1(4(1)))&2&0.14&1.8&1&(312)[+\,+\,+]\\ \hline
(1(4(4)))&2&0.999&0.4975&0.9&(213)[--\,+\,+]\\ \hline
(1(5(4)))&3&0.997&2.15&1.3&(312)[--\,+\,+]\\ \hline
(1(5(5)))&0.2&0.8&1&1&(213)[+\,+\,+]\\ \hline
(2(1(3)))&0.1&0.1&0.2&1.7&(213)[+\,+\,--]\\ \hline
(2(1(4)))&1.4&0.62&4.6&2&(312)[+\,+\,+]\\ \hline
(2(2(2)))&1.7&0.08&0.3&2.7&(123)[--\,+\,--]\\ \hline
(2(2(3)))&0.9&0.55&2.65&4&(213)[+\,--\,+]\\ \hline
\end{tabular}\]}
\begin{rem} To each algebraic code correspond several geometrical codes
and vice versa. The only exception here is: the geometrical code $(1,2,3)$
corresponds only to the algebraic code $(123)[+\,+\,+]$.
\end{rem}

\section{Supplement 1}
\pn Here we present a sketch of the proof of Theorem 3.1. \pmn Let
$p_1=x^2+a^2$, $p_2=(x-b)^2+c^2$, $p_3=(x-1)^2+d^2$ and $p=\int
p_1p_2p_3dx$. We assume that $a,c,d>0$ and $0<b<1$. Let
$x_1=a\,i$, $x_2=b+c\,i$ and $x_3=1+d\,i$ --- be critical points
and $p(x_1)=r_1+i\,q_1$, $p(x_2)=r_2+i\,q_2$ and
$p(x_3)=r_3+i\,q_3$
--- be corresponding critical values. Let
$$p(x)=Q_j(x)\cdot p_j^2(x)+W_j(x),\,\,j=1,2,3. \eqno(1)$$ Then
$$p'(x)=Q_j'(x)\cdot p_j^2(x)+2Q_j(x)\cdot p_j(x)\cdot p_j'(x)+
W_j'(x). \eqno(2)$$ As $p'(x_j)=0$ and $p_j(x_j)=0$, then $W_j'(x_j)=0$
(and $W_j'(\overline{x_j})=0$). As degree $W_j'$=2, then
$W_j'(x)=k_j\cdot p_j(x)$. As $p(x_j)=W_j(x_j)$, then it is easy to check that
$q_1=\frac 23\cdot a^3k_1$, $q_2=\frac 23\cdot c^3k_2$,
$q_3=\frac 23\cdot d^3k_3$. Thus, signs of $q_1,q_2,q_3$ are the same, as
signs $k_1,k_2,k_3$, respectively. We have
$$175\cdot k_2=8c^4+140c^2b(1-b)+7\cdot L\cdot R, \eqno(3)$$
where $L=5a^2-c^2+5b^2$ and $R=5d^2-c^2+5b(1-b)$. Thus, $k_2<0$
only if $LR<0$. \pmn \underline{The case $k_2<0,L<0,R>0$.} Let
$t=-L=c^2-5a^2-5b^2>0$. Thus, $c^2=t+5a^2+5b^2$. We replace $c^2$
in $k_1$ and $k_2$ by $t+5a^2+5b^2$. We have
$$175\,t\,k_1=175(24a^2+30b^2+5t)\cdot(-k_2)+V_1, \eqno(4)$$
were  $V_1=\sum_{j=0}^6 \binom{6}{j}\alpha_jb^j(1-b)^{6-j}$ and
{\small $$\begin{array}{l}
\alpha_0= 105(4a^2+t)^3\\
\alpha_1= 35(4a^2+t)^2(12a^2+3t+7)\\
\alpha_2= 21(4a^2+t)(80a^4+40a^2t+92a^2+5t^2+24t)\\
\alpha_3= 6720a^6+5040a^4t+11424a^4+1260a^2t^2+5964a^2t+2646a^2+
105t^3+777t^2+686t\\
\alpha_4= 6720a^6+5040a^4t+15008a^4+1260a^2t^2+8008a^2t+8624a^2+
105t^3+1064t^2+7028/3t\\
\alpha_5= 6720a^6+5040a^4t+18480a^4+1260a^2t^2+10080a^2t+
16660a^2+105t^3+1365t^2+14560/3t+4900\\
\alpha_6= 6720a^6+5040a^4t+21840a^4+1260a^2t^2+12180a^2t+
23520a^2+105t^3+1680t^2+7700t+8400 \end{array}$$}
Thus, $k_1>0$. To prove the second part of the statement, i.e.
$r_1>r_2 \Leftrightarrow {\rm Re}(p(a\,i))>{\rm Re}(p(b+c\,i))$,
we note that
$$470\cdot t(r_1-r_2)=245b(240a^4+540a^2b^2+120a^2t+324b^4+140b^2t+15t^2)
\cdot(-k_2)+\sum_{j=0}^9 \binom{9}{j}\beta_j b^j(1-b)^{9-j},$$ where
{\scriptsize
$$\begin{array}{l}
\beta_0=245t(4a^2+t)^3\\ \\
\beta_1=\frac{70}{9}\,(4a^2+t)^2(60a^4+122a^2t+29t^2)\\ \\
\beta_2=\frac{35}{9}\,(4a^2+t)(960a^6+1184a^4t+420a^4
+448a^2t^2+252a^2t +53t^3+42t^2)\\ \\ \beta_3=\frac{10}{3}\,(
6720a^8+7616a^6t+6060a^6 +3444a^4t^2+4794a^4t +714a^2t^3 +
1347a^2t^2 +56t^4+132t^3)\\
\\ \beta_4=\frac{5}{9}\,(53760a^8+51520a^6t+74880a^6 +20496a^4t^2
+55080a^4t+10920a^4 + 3948a^2t^3 +14688a^2t^2 +5544a^2t
+301t^4+1404t^3+756t^2)\\ \\
\beta_5=\frac{10}{9}\,(33600a^8+28672a^6t+64200a^6 +10164a^4t^2
+43860a^4t+25128a^4 + 1806a^2t^3 +11010a^2t^2 +12204a^2t
+133t^4+1020t^3+1620t^2)\\ \\
\beta_6=\frac{5}{3}\,(26880a^8+21056a^6t+66000a^6 +6720a^4t^2
+41820a^4t+45912a^4 + 1092a^2t^3 +9780a^2t^2 +21096a^2t+6048a^2+\\
\hspace{2cm}+77t^4+870t^3 +2700t^2+1512t)\\ \\
\beta_7=\frac{70}{9}\,(6720a^8+4928a^6t+20340a^6 +1428a^4t^2
+11934a^4t+ 20784a^4 +210a^2t^3 +2565a^2t^2 +8892a^2t+7236a^2+\\
\hspace{2cm}+14t^4+216t^3 +1080t^2+1728t)\\ \\
\beta_8=\frac{35}{9}\,(15360a^8+10688a^6t+55680a^6 +2832a^4t^2
+30192a^4t+74448a^4 + 372a^2t^3 +5856a^2t^2 +28944a^2t+43200a^2+\\
\hspace{2cm}+23t^4+456t^3 +3240t^2+9504t+9072)\\ \\
\beta_9=70\,(960a^8+640a^6t+4080a^6 + 156a^4t^2
+2040a^4t+6576a^4 +18a^2t^3 + 348a^2t^2 +2232a^2t+4752a^2+\\
\hspace{2cm} +t^4+24t^3 +216t^2 +864t+1296) \end{array}$$} \pmn
\underline{The case $k_2<0,L>0,R<0$} can be considered
analogously.

\section{Supplement 2}
\pn Below is (computed by us) the table of realized (marked by
"\checkmark") algebraic codes for the 9-th degree polynomials.
Here "p"{} is for positive and "n"{} --- for negative. {\tiny
\[\begin{tabular}{|c|c|c|c|c|c|c|c|c|c|c|c|c|c|c|c|c|}
\hline &pppp&pppn&ppnp&ppnn&pnpp&pnpn&pnnp&pnnn&
nppp&nppn&npnp&npnn&nnpp&nnpn&nnnp&nnnn\\
\hline

(1234) & \checkmark & \checkmark & \checkmark & \checkmark &
\checkmark & \checkmark & 0 & 0 & \checkmark & \checkmark &
\checkmark & \checkmark & \checkmark & \checkmark & 0 & 0 \\
\hline

(1243) & \checkmark & \checkmark & \checkmark & \checkmark &
\checkmark & \checkmark & \checkmark & 0 & \checkmark & \checkmark
& \checkmark & \checkmark & \checkmark & \checkmark & \checkmark & 0 \\
\hline

(1324) & \checkmark & \checkmark & \checkmark & \checkmark &
\checkmark & \checkmark & 0 & 0 & \checkmark & \checkmark &
\checkmark & \checkmark & \checkmark & \checkmark & 0 & 0 \\
\hline

(1342) & \checkmark & \checkmark & \checkmark & \checkmark &
\checkmark & \checkmark & 0 & 0 & \checkmark & \checkmark &
\checkmark & \checkmark & \checkmark & \checkmark & \checkmark & 0 \\
\hline

(1423) & \checkmark & \checkmark & \checkmark & \checkmark &
\checkmark & \checkmark & \checkmark & 0 & \checkmark & \checkmark
& \checkmark & \checkmark & \checkmark & \checkmark & \checkmark & 0 \\
\hline

(1432) & \checkmark & \checkmark & \checkmark & \checkmark &
\checkmark & \checkmark & \checkmark & 0 & \checkmark & \checkmark
& \checkmark & \checkmark & \checkmark & \checkmark & \checkmark & 0 \\
\hline

(2134) & \checkmark & \checkmark & \checkmark & \checkmark &
\checkmark & \checkmark & \checkmark & \checkmark & \checkmark &
\checkmark & \checkmark & \checkmark & \checkmark & \checkmark & 0 & 0 \\
\hline

(2143) & \checkmark & \checkmark & \checkmark & \checkmark &
\checkmark & \checkmark & \checkmark & \checkmark & \checkmark &
\checkmark & \checkmark & \checkmark & \checkmark & \checkmark & \checkmark & 0 \\
\hline

(2314) & \checkmark & \checkmark & \checkmark & \checkmark &
\checkmark & \checkmark & \checkmark & \checkmark & \checkmark &
\checkmark & \checkmark & \checkmark & \checkmark & \checkmark & 0 & 0 \\
\hline

(2341) & \checkmark & \checkmark & \checkmark & \checkmark &
\checkmark & \checkmark & \checkmark & \checkmark & \checkmark &
\checkmark & \checkmark & \checkmark & \checkmark & \checkmark & 0 & 0 \\
\hline

(2413) & \checkmark & \checkmark & \checkmark & \checkmark &
\checkmark & \checkmark & \checkmark & \checkmark & \checkmark &
\checkmark & \checkmark & \checkmark
& \checkmark & \checkmark & \checkmark & 0 \\
\hline

(2431) & \checkmark & \checkmark & \checkmark & \checkmark &
\checkmark & \checkmark & \checkmark & \checkmark & \checkmark &
\checkmark & \checkmark & \checkmark
& \checkmark & \checkmark & \checkmark & 0 \\
\hline

(3124) & \checkmark & \checkmark & \checkmark & \checkmark &
\checkmark & \checkmark & 0 & \checkmark & \checkmark & \checkmark
& \checkmark & \checkmark & \checkmark & \checkmark & 0 & 0 \\
\hline

(3142) & \checkmark & \checkmark & \checkmark & \checkmark &
\checkmark & \checkmark & 0 & \checkmark & \checkmark & \checkmark
& \checkmark & \checkmark & \checkmark & \checkmark & \checkmark & 0 \\
\hline

(3214) & \checkmark & \checkmark & \checkmark & \checkmark &
\checkmark & \checkmark & \checkmark & \checkmark & \checkmark &
\checkmark & \checkmark & \checkmark & \checkmark & \checkmark & 0 & 0 \\
\hline

(3241) & \checkmark & \checkmark & \checkmark & \checkmark &
\checkmark & \checkmark & \checkmark & \checkmark & \checkmark &
\checkmark & \checkmark & \checkmark & \checkmark & \checkmark & 0 & 0 \\
\hline

(3412) & \checkmark & \checkmark & \checkmark & \checkmark &
\checkmark & \checkmark & \checkmark & \checkmark & \checkmark &
\checkmark & \checkmark & \checkmark & \checkmark & \checkmark & \checkmark & 0 \\
\hline

(3421) & \checkmark & \checkmark & \checkmark & \checkmark &
\checkmark & \checkmark & \checkmark & \checkmark & \checkmark &
\checkmark & \checkmark & \checkmark & \checkmark & \checkmark & \checkmark & 0 \\
\hline

(4123) & \checkmark & \checkmark & \checkmark & \checkmark &
\checkmark & \checkmark & \checkmark & 0 & \checkmark & \checkmark
& \checkmark & \checkmark & \checkmark & \checkmark & \checkmark & 0 \\
\hline

(4132) & \checkmark & \checkmark & \checkmark & \checkmark &
\checkmark & \checkmark & \checkmark & 0 & \checkmark & \checkmark
& \checkmark & \checkmark & \checkmark & \checkmark & \checkmark & 0 \\
\hline

(4213) & \checkmark & \checkmark & \checkmark & \checkmark &
\checkmark & \checkmark & \checkmark & \checkmark & \checkmark &
\checkmark & \checkmark & \checkmark
& \checkmark & \checkmark & \checkmark & 0 \\
\hline

(4231) & \checkmark & \checkmark & \checkmark & \checkmark &
\checkmark & \checkmark & \checkmark & \checkmark & \checkmark &
\checkmark & \checkmark & \checkmark & \checkmark & \checkmark & \checkmark & 0 \\
\hline

(4312) & \checkmark & \checkmark & \checkmark & \checkmark &
\checkmark & \checkmark & \checkmark & \checkmark & \checkmark &
\checkmark & \checkmark & \checkmark & \checkmark & \checkmark & \checkmark & 0 \\
\hline

(4321) & \checkmark & \checkmark & \checkmark & \checkmark &
\checkmark & \checkmark & \checkmark & \checkmark & \checkmark &
\checkmark & \checkmark & \checkmark & \checkmark & \checkmark & \checkmark & 0 \\
\hline \end{tabular}\]}  \pn It must be noted that zeroes in this
table are only in columns marked by two middle "n".

\vspace{5mm}
\end{document}